\documentclass[12pt,a4paper]{amsart}
\usepackage{amsmath,amsfonts,amssymb,amsthm,amscd}
\newtheorem{thm}{Theorem}
\newtheorem{lem}[thm]{Lemma}
\newtheorem{defn}[thm]{Definition}
\newtheorem{prop}[thm]{Proposition}

\begin{document}

\title{Sasaki--Weyl connections on CR manifolds}
\date{}
\author{Liana David}
\address{Mathematics
Institute of the Romanian Academy ``Simion Stoilow'', Calea
Grivi\c{t}ei nr. 21, Bucharest, Romania}
\email{liana.david@imar.ro}
\begin{abstract}
  We introduce and study the notion of Sasaki--Weyl manifold, which is
  a natural generalization of the notion of Sasaki manifold. We
  construct a reduction of Sasaki--Weyl manifolds and we show that it
  commutes with several reductions already existing in the literature.
\end{abstract}
\maketitle

\section{Introduction}

The reduction construction is well-known in many areas of
differential geometry. It has been initially introduced in
symplectic geometry by Marsden and Weinstein and then extended in
the context of K\"{a}hler, hyper-K\"{a}hler, CR and Sasaki
geometry. Reductions in conformal geometry have also been studied:
a reduction for symplectic conformal manifolds appears in
\cite{haller} and a reduction for locally conformal K\"{a}hler
manifolds has been constructed in \cite{gaud}.

The aims of this paper are threefold: (1) to introduce and study
the notion of Sasaki--Weyl manifold, which generalizes, in the case
when the Weyl connection is closed, the notion of Sasaki manifold,
in the same way as the locally conformal K\"{a}hler manifolds
generalize the K\"{a}hler manifolds; (2) to construct and
investigate the properties of a reduction for Sasaki--Weyl
manifolds, which completes the picture of reductions already
existing in the literature; (3) to illustrate our theory with
examples and indicate further directions of possible research.

Let $(M,c)$ be a conformal oriented manifold. H.Weyl noticed (see
\cite{weyl}) that there is a one-to-one correspondence between
connections on the density line bundle $L$ of $M$ (usually called
Weyl connections) and torsion-free connections on $M$ which
preserve the conformal structure $c$. The density line bundle $L$
is oriented and positive sections of $L$ correspond to metrics in
the conformal class $c$. If the Weyl connection $D$ is exact (i.e.
there is a global non-vanishing $D$-parallel section $\mu$ of $L$)
then the corresponding connection on $M$ is the Levi-Civita
connection of the corresponding Riemannian metric $\mu^{-2}c$. 
Consider now an oriented (strongly pseudo-convex) CR
manifold $(M,H,I)$ and $L:=TM/H$ its co-contact bundle, which
inherits an orientation from the orientation of $M$ and $H$.
Positive sections of $L$ correspond to contact forms of $M$.
Tanaka associated (see \cite{tanaka}), in a natural way, to every
contact form of $H$ a connection on $M$, the so called Tanaka
connection, which preserves the contact form and the CR structure
of $M$ and satisfies other natural conditions. Tanaka's result has
been generalized in \cite{liana}, where to every connection $D$ on
$L$ it was associated, in a natural way, a connection on $M$ which
coincides with the Tanaka connection of a contact form in the case
when $D$ is exact. Therefore, there is an obvious similarity
between conformal and CR geometry, Weyl's theorem from
conformal geometry corresponding to the generalized Tanaka's
theorem from CR geometry. By analogy, connections on the
co-contact bundle $L$ of a CR manifold $(M,H,I)$ were called in
\cite{liana} Weyl connections.

In this paper we will push the similarity between conformal and CR
geometry further in the realm of K\"{a}hler and Sasaki geometry. A
Sasaki--Weyl manifold is a CR oriented manifold $(M,H,I)$ together
with a Weyl connection $D$ whose Reeb vector field $\psi^{D}$
satisfies the Sasaki--Weyl condition $L_{\psi^{D}}(I)=0.$ 
When $D$ is exact the Sasaki--Weyl condition is precisely the 
Sasaki condition for the contact form associated to a positive
$D$-parallel section of the co-contact bundle $L$ of $(H,I).$ When $D$ is
closed, the Sasaki--Weyl condition translates to the Sasaki
condition on the universal cover of $M$ with the CR
structure induced by $(H,I)$. After briefly recalling some facts we need 
from CR
geometry, we show in Section \ref{structcomplexa} that every Weyl
connection $D$ on an oriented CR manifold $(M,H,I)$ determines, in
a natural way, an almost hermitian conformal structure on the cone
manifold $L^{+}$ of oriented elements of the co-contact bundle $L$
of $M$. It turns out that this almost hermitian structure is
locally conformal K\"{a}hler when $D$ is Sasaki--Weyl with
curvature $F^{D}$ satisfying an additional condition.  
In Section \ref{reducere} we perform a reduction of Sasaki--Weyl
manifolds and then we study its properties. While closed
Sasaki--Weyl connections remain closed through the reduction
process, this is not true for exact Sasaki--Weyl connections.
Proposition \ref{uni} shows that in fact any Sasaki--Weyl manifold
whose Weyl connection is closed can be written as a reduction of a
Sasaki--Weyl manifold whose Weyl connection is exact. We show the
compatibility of the Sasaki--Weyl reduction with the operation of
taking universal covers and in Section \ref{comutativitate} we
show its commutativity with the locally conformal K\"{a}hler
reduction of Biquard-Gauduchon of the corresponding cone
manifolds. We end the paper with examples of closed and non-closed
Sasaki--Weyl connections, and an explicit illustration of the
Sasaki--Weyl reduction. We hope that our theory (in particular the
cone construction we propose in the paper) will provide examples
of Vaisman manifolds. Further investigation in this direction is
needed.

\paragraph{Acknowledgements:} I warmly thank Toby Bailey for
useful discussions related to this paper, David Calderbank for
suggesting me to study Weyl connections on CR manifolds and Paul
Gauduchon, for encouraging me to work on CR geometry and for
kindly giving me his unpublished notes on the locally conformal
K\"{a}hler reduction. I thank Liviu Ornea for drawing my attention
to references \cite{loose1}, \cite{loose2} and EPSRC for financial
support.

\section{Sasaki--Weyl manifolds}

We recall that on an oriented (strongly pseudo-convex) CR manifold
$(M,H,I)$, $H$ is a codimension one oriented subbundle of $TM$,
$I:H\rightarrow H$ is an endomorphism of $H$ which satisfies
$I^{2}=-\mathrm{Id}$ (here and elsewhere \lq\lq$\mathrm{Id}$" denotes
the identity endomorphism) and the following integrability
condition holds: if $X$ and $Y$ are sections of $H$ then $[IX,IY]-[X,Y]$
is a section of $H$ as well and the relation
$$
[IX,IY]-[X,Y]=I([IX,Y]+[X,IY])
$$
holds. The co-contact bundle $L:=TM/H$ is oriented and a positive
section of its dual $L^{-1}$ is a contact form. To every contact
form $\theta$ (viewed as a $1$-form on $M$ with kernel $H$) there
is associated a Reeb vector field $T$, determined by the
conditions $\theta (T)=1$ and $i_{T}d\theta =0$, and a metric
$g_{0}$ on $H$ defined by $g_{0}:=\frac{1}{2}d\theta (\cdot
,I\cdot )$, which is independent of the choice of $\theta$ and
must be positive definite due to the strongly pseudo-convexity
condition. If $X$ is a section of $H$, so is $L_{T}(X)$ and
$L_{T}(I)$, defined by $L_{T}(I)(X):=L_{T}(IX)-IL_{T}(X)$, is an
endomorphism of $H$. If $L_{T}(I)=0$ we say that $T$ (and
$\theta$) is the Reeb
vector field (contact form respectively) of a Sasaki structure.\\

We now briefly recall some basic facts from the theory of CR-Weyl
manifolds, as developed in \cite{liana}. A CR-Weyl manifold is an
oriented CR manifold $(M,H,I)$ together with a connection $D$
(called a Weyl connection) on its co-contact bundle $L.$ The
curvature $F^{D}$ of $D$ is called the Faraday curvature. The
connection $D$ is closed if $F^{D}=0$ and exact if there is a
non-vanishing $D$-parallel section of $L$. Let $\eta
:TM\rightarrow L$ be the natural projection (sometimes called the
universal contact form). The Weyl connection $D$ determines a Reeb
vector field $\psi^{D}:L\rightarrow TM$, which 
is a bundle homomorphism (or a section of the bundle $L^{-1}\otimes TM$)
uniquely determined by the two conditions
$i_{\psi^{D}}d^{D}\eta =0$ and $\eta\circ\psi^{D}=\mathrm{Id}$.
The terminology \lq\lq Reeb vector field" for $\psi^{D}$ is justified by
the following fact: if $D$ is exact and $\mu$ is a $D$-parallel
positive section of $L$, then $\psi^{D}(\mu )$, which is a genuine
vector field on $M$, is the Reeb vector field of the contact form
$\theta :=\mu^{-1}\eta$ (for a non-vanishing section $\mu$ of $L$
we denote by $\mu^{-1}$ the section of $L^{-1}$ 
induced by $\mu$, i.e. the natural contraction between $\mu$ and $\mu^{-1}$ 
is the function identically one on $M$).

Associated to the Weyl connection $D$ there is also a Weyl-Lie
derivative $L_{\psi^{D}}$, which on vector fields $X$ on $M$ is
defined by
$$
L_{\psi^{D}}(X)=\mu^{-1}\left( L_{\psi^{D}(\mu )}(X)+D_{X}(\mu
)\psi^{D}\right)
$$
and which is independent of the chosen (non-vanishing) section $\mu$ of
$L$. If $X$ is a section of $H$, then $L_{\psi^{D}}(X)$ is a
section of $L^{-1}\otimes H$. Due to this observation, we can
define $L_{\psi^{D}}(I)$, an endomorphism of $H$ with values in
the bundle $L^{-1}$, by the formula
$L_{\psi^{D}}(I)(X):=L_{\psi^{D}}(IX)-IL_{\psi^{D}}(X)$.

\begin{defn}
A CR-Weyl manifold $(M,H,I,D)$ whose Weyl connection $D$ satisfies
$L_{\psi^{D}}(I)=0$ is called  Sasaki--Weyl.
\end{defn}

\textbf{Remark:} We explain now the analogy between Sasaki and
Sasaki--Weyl manifolds with closed Weyl connection on one hand, and
K\"{a}hler and locally conformal K\"{a}hler manifolds on the other
hand. Recall that an (oriented) locally conformal K\"{a}hler
manifold (of real dimension greater or equal to $4$) is an almost
hermitian conformal manifold $(M,c,J)$ with a connection $D^{c}$
(called the canonical connection) on the density line bundle $L$
of $M$ which satisfies two properties: $D^{c}$ is closed (i.e. its
curvature is zero); the corresponding connection on $M$
(obtained by Weyl's theorem) preserves $J$. This implies, in
particular, that $J$ is integrable. Conversely, if $J$ is
integrable and $\Omega$ is the K\"{a}hler form of an arbitrary
metric $\mu^{-2}c$ in the conformal class $c$, then $(M,c,J)$ is
locally conformal K\"{a}hler if there is a (unique) closed
$1$-form $\alpha$ on $M$ such that $d\Omega =\Omega\wedge\alpha$.
Moreover, the canonical connection $D^{c}$ of $(M,c,J)$ satisfies
$D_{X}^{c}(\mu )=-\frac{1}{2}\alpha (X)\mu$, for every vector
field $X$ on $M$. Every local, positive, $D^{c}$-parallel section
of $L$ generates a local K\"{a}hler metric in the conformal class
$c$. On the other hand, a Sasaki--Weyl manifold $(M,H,I,D)$ whose 
Weyl connection $D$ is closed enjoys a similar property: every local, 
positive, $D$-parallel section $\mu$ of the co-contact bundle of $(M,H,I)$ 
determines a local contact form $\mu^{-1}\eta$ of a Sasaki structure.\\

In this paper we shall use the following notations and
conventions:

\paragraph{Notation:}
Unless otherwise specified, all our manifolds will be connected
and oriented. We shall denote by ${\mathcal X}(M)$ the space of
vector fields on a manifold $M$. For a CR-Weyl manifold
$(M,H,I,D)$ we shall use the following notations: $\pi
:L\rightarrow M$ for its co-contact bundle; 
$\Lambda$ for the image of the Reeb vector field $\psi^{D}:L\rightarrow
TM$, which is a rank one subbundle of $TM$ complementary to the contact bundle
$H$; $\bar{D}$ for the
connection on the bundle $\pi^{*}L\rightarrow L$ induced by the
Weyl connection $D$; $F^{D,-}$ for the $I$-anti-invariant part of
$F^{D}\vert_{H\times H}$; $\tilde{X}$ for the $D$-horizontal lift
to $L$ of $X\in{\mathcal X}(M).$

\section{The cone of a Sasaki--Weyl manifold}\label{structcomplexa}

In this Section we fix a CR-Weyl manifold $(M,H,I,D)$.

\begin{prop}\label{complexa}
Let $y\in L.$ The connection $D$ determines a decomposition
\begin{equation}\label{desc} 
T_{y}L=H_{\pi (y )}\oplus\Lambda_{\pi (y)}\oplus  L_{\pi (y )},
\end{equation}
and an almost complex structure $J$ on the manifold $L$, which, by means 
of the decomposition (\ref{desc}), is written as
$$
J|_{H}=I, \quad J(s):=\psi^{D}(s),
$$
where $s\in L_{\pi (y)}$ is non-zero. The almost complex
structure $J$ is integrable if and only if $D$ is Sasaki--Weyl,
$F^{D}|_{H\times H}$ is $I$-invariant and $i_{\psi^{D}}F^{D}=0$.
\end{prop}

\begin{proof}
Recall that if a section $s$ of $L$ is viewed as a vertical vector
field on $L$ and $X,Y\in{\mathcal X}(M)$, then $[\tilde{X},s
]=D_{X}(s)$ and $v^{D}([\tilde{X},\tilde{Y}]_{y})=-F^{D}_{X,Y}y$,
where "$v^{D}$" denotes the projection onto the last factor in
the decomposition (\ref{desc}). Let $X$ and $Y$ be two sections of $H$, viewed
as horizontal vector fields on $L$. 
The integrability tensor $N^{J}$ of $J$, applied to the pair $(X,Y)$, can be
calculated as follows: 
\begin{align*}
N^{J}(X,Y)_{y}&=[J\tilde{X},J\tilde{Y}]_{y}-J\left( [J\tilde{X},\tilde{Y}]_{y}
+[\tilde{X},J\tilde{Y}]_{y}\right) -[\tilde{X},\tilde{Y}]_{y}\\
&=\widetilde{[IX,IY]}_{y}-F^{D}_{IX,IY}y-J\left( \widetilde{[IX,Y]}_{y}+
\widetilde{[X,IY]}_{y}-F^{D}_{IX,Y}y-F^{D}_{X,IY}y\right)\\
&-\widetilde{[X,Y]}_{y}+F^{D}_{X,Y}y\\
&=-F^{D}_{IX,IY}y+J\left( F^{D}_{IX,Y}y+F^{D}_{X,IY}y\right) +F^{D}_{X,Y}y\\
&=-2F^{D,I-}_{X,Y}y
+2F^{D,I-}_{IX,Y}\widetilde{\psi^{D}(y)},
\end{align*}
where we have used the integrability of the complex structure $I$ of the 
bundle $H$. A similar calculation shows that
$$
N^{J}(\tilde{X},\widetilde{\psi^{D}(s)})_{y}= F^{D}_{X,\psi^{D}(s)}y
+F^{D}_{IX,\psi^{D}(s)}\widetilde{\psi^{D}(y)}-
\widetilde{sL_{\psi^{D}}(I)(IX)},
$$
where $s$ is a local smooth section of $L$ non-vanishing
at the point $\pi (y).$ The conclusion follows.
\end{proof}

\begin{lem}\label{potential}
Let $\sigma$ be the natural section of $\pi^{*}L\rightarrow L$ and
$\omega$ the $\pi^{*}L^{2}$-valued $2$-form on $L$ defined by
$\omega :=\frac{1}{2}d^{\bar{D}}(\sigma\pi^{*}\eta ).$ Then
$\omega$ has $J$-potential $\sigma^{2} .$ Equivalently, $\omega
=\frac{1}{4}d^{\bar{D}}\left( J\bar{D}\sigma^{2}\right) .$
\end{lem}

\begin{proof}
Let $s_{0}$ be a positive section of $L$, $\gamma
:=s_{0}^{-1}Ds_{0}$, $\theta_{0}:=s_{0}^{-1}\eta$ and $T_{0}:=
\psi^{D}(s_{0}).$  The section
$s_{0}$ determines a trivialization $L\cong M\times\mathbb{R}$ of
$L$. In the trivialization of $\pi^{*}L^{2}$ defined by
$\pi^{*}s_{0}^{2}$, $\sigma^{2}(x,t)=(x,t,t^{2})$ and
$\bar{D}\sigma^{2}=(dt^{2}+2t^{2}\pi^{*}\gamma )\pi^{*}s_{0}^{2}.$
In order to calculate $J\bar{D}(\sigma^{2})$, we need to determine
$J\pi^{*}\gamma$ and $d^{J}t^{2}.$ Note that, for $X\in H_{x}$,
$J\tilde{X}_{(x,t)}=IX-t\gamma (IX)\frac{\partial}{\partial t}$
and $J\widetilde{T}
=-s_{0}=-\frac{\partial}{\partial t}$, which easily imply
that
\begin{align*}
JX&=IX+t\gamma (X)T_{0}-[t^{2}\gamma (X)\gamma 
(T_{0})+t\gamma (IX)]\frac{\partial}{\partial t}\\
J\left( \frac{\partial}{\partial t}\right)&=
T_{0}-t\gamma (T_{0})\frac{\partial}{\partial t},
\end{align*}
for every $X\in H$. Then  
$$
J\pi^{*}\gamma =\pi^{*}J_{0}\gamma -t\gamma (T_{0})
\pi^{*}\gamma -\gamma (T_{0})dt,
$$
where $J_{0}\in\mathrm{End}(TM)$ extends $I\in\mathrm{End}(H)$
being zero on $T_{0}$, and
$$
d^{J}t=-t\pi^{*}(J_{0}\gamma )+t^{2}\gamma
(T_{0})\pi^{*}\gamma +\pi^{*}(\theta_{0})+t\gamma
(T_{0})dt.
$$
It follows that $J\left(dt^{2}+2t^{2}\pi^{*}\gamma\right)
=2t\pi^{*}\theta_{0}$ or equivalently
$J\bar{D}(\sigma^{2})=2\sigma\pi^{*}\eta .$
\end{proof}

On the total space $L$ of the co-contact bundle of $M$ we have an
almost complex structure $J$, provided by Proposition
\ref{complexa} and a conformal $2$-form $\omega$ provided by Lemma
\ref{potential}. Their properties are expressed by the following
Lemma, whose proof is an easy calculation.

\begin{lem}\label{metrica}
The $2$-form $\omega$ is $J$-invariant. In the trivialization of
$L$ determined by the positive section $s_{0}$, the bilinear form
$g:=\pi^{*}(s_{0}^{-2})\omega (\cdot ,J\cdot )$ has the following
expression:
$g(\tilde{X},\tilde{Y})=\frac{t}{2}d(s_{0}^{-1}\eta)(X,IY)$, for
$X,Y\in H$;
$g(\widetilde{\psi^{D}(s_{0})},\widetilde{\psi^{D}(s_{0})})=
g\left(\frac{\partial}{\partial t},\frac{\partial}{\partial
t}\right)=\frac{1}{2}$;
$g(\tilde{X},\widetilde{\psi^{D}(s_{0})})=0$, for $X\in H.$ In
particular, $g$ is positive definite (with values in the oriented
line bundle $\pi^{*}L^{2}$) on the subset $L^{+}$ of positive
oriented elements of $L$. Moreover, for every $y\in L^{+}$ the
decomposition
$$
T_{y}L^{+}=H_{\pi (y )}\oplus\Lambda_{\pi (y)}
\oplus L_{\pi (y)}
$$
is $g$-orthogonal.
\end{lem}

\begin{defn}\label{cone}
The cone of the CR-Weyl manifold $(M,H,I,D)$ is the manifold
$L^{+}$ together with the almost complex structure $J$ and the
$\pi^{*}L^{2}$-valued $2$-form $\omega$.
\end{defn}

\begin{thm}\label{lcK}
Let $(L^{+},\omega ,J)$ be the cone of the oriented CR-Weyl
manifold $(M,H,I,D)$. Then $(L^{+},\omega ,J)$ is locally
conformal K\"{a}hler (l.c.K.) if and only if $D$ is Sasaki--Weyl
and $F^{D}=kd^{D}\eta$, for a section $k$ of $L^{-1}$ which
satisfies $D(k)+k^{2}\eta =0.$ Moreover, on the subset of $L^{+}$
where $\pi^{*}(k)$ is non-vanishing, $(L^{+},\omega ,J)$ is
globally conformal K\"{a}hler with K\"{a}hler form
$\pi^{*}(k^{2})\omega$.
\end{thm}

\begin{proof}
With the notation of Lemma \ref{potential}, let $\Omega
:=\pi^{*}(s_{0}^{-2})\omega =\frac{1}{2}d(t\pi^{*}\theta_{0})+t
\pi^{*}(\gamma\wedge\theta_{0})$. The exterior derivative of
$\Omega$ satisfies $d\Omega =t\pi^{*}(F^{D}\wedge\theta_{0})
-2\Omega\wedge\pi^{*}\gamma .$ It follows that $(L^{+},\omega ,J)$
is locally conformal K\"{a}hler if $J$ is integrable and if there
is a $1$-form $\alpha :=k_{0}dt+k_{1}$ on $L^{+}$, where $k_{0}\in
C^{\infty}(M\times\mathbb{R}^{>0})$ and $k_{1}$ is a $1$-form on
$M$ parametrized by $t\in \mathbb{R}^{>0}$, such that the
following two relations:
$\pi^{*}(F^{D}\wedge\theta_{0})=\Omega\wedge\alpha$  and
$d(t\alpha -2\pi^{*}\gamma )=0$ hold. The first of these relations
is equivalent to the system
\begin{align*}
\frac{k_{0}t}{2}\pi^{*}(d\theta_{0})+tk_{0}\pi^{*}(\gamma\wedge\theta_{0})
-\frac{1}{2}k_{1}\wedge\pi^{*}\theta_{0}&=0\\
\pi^{*}(F^{D}\wedge\theta_{0})-[\frac{t}{2}\pi^{*}(d\theta_{0})+t\pi^{*}(\gamma\wedge\theta_{0})]\wedge
k_{1}&=0.
\end{align*}
Since $\theta_{0}$ annihilates the vectors in $H$ and
$d\theta_{0}$ is non-degenerate when restricted to $H$, the first
equation implies that $k_{0}=0$ and
$k_{1}\wedge\pi^{*}\theta_{0}=0.$ This means that $\alpha =k_{1}=
\lambda\pi^{*}\theta_{0}$ for a function $\lambda\in
C^{\infty}(M\times\mathbb{R}^{>0})$ which implies that
$\pi^{*}(F^{D}\wedge\theta_{0})
=\frac{t\lambda}{2}\pi^{*}(d\theta_{0}\wedge\theta_{0})$, or
$F^{D}|_{H\times H}=\frac{t\lambda s_{0}^{-1}}{2}d^{D}\eta
|_{H\times H}.$ Since $J$ is integrable, $i_{\psi^{D}}F^{D}=0$
(see Proposition \ref{complexa}) and we obtain that
$F^{D}=kd^{D}\eta $, for $k:=\frac{t\lambda s_{0}^{-1}}{2}$ which
must be independent of $t$. On the other hand, the relation
$d\left({t\alpha}-2\pi^{*}\gamma\right) =0$ is equivalent with
$F^{D}=d(k\eta )$. Since $F^{D}=d(k\eta )=kd^{D}\eta$, the section
$k$ of $L^{-1}$ must satisfy $D(k)\wedge\eta =0$. Also,
$d(kd^{D}\eta )=\left( D(k)+k^{2}\eta\right)\wedge d^{D}\eta =0$,
which in turn implies, since $i_{\psi^{D}}d^{D}\eta =0$, that
$D_{\psi^{D}}(k)+k^{2}=0$. The first claim of the Theorem follows.
For the second claim, we notice that the canonical Weyl connection
$D^{c}$ of $(L^{+},\omega ,J)$ is $D^{c}=\pi^{*}(D-k\eta )$ and
satisfies $D^{c}(\pi^{*}k^{-1})=0.$
\end{proof}

\section{Reductions of Sasaki--Weyl manifolds}\label{reducere}

We recall that a Lie group $G$ acts by CR automorphisms on a CR
manifold $(M,H,I)$ if it preserves $H$ and $I$:
$g_{*}(H_{x})\subset H_{gx}$ and $g_{*}IX=Ig_{*}X$ for every $g\in
G$, $x\in M$ and $X\in H_{x}$. Since $G$ preserves $H$, there is
an induced action of $G$ on the co-contact bundle $L$ of $M$ and
on the space $\Gamma (M,L)$ of smooth sections of $L$. We shall
usually denote by $g\cdot y$ the action of  $g\in G$ on $y\in L$.
The action of $G$ on $\Gamma (M,L)$ is defined by $(g\cdot
s)(x)=g\cdot s(g^{-1}x)$ for every $g\in G$, $s\in\Gamma (M,L)$
and $x\in M.$ Any vector field $\xi^{a}$ on $M$, generated by an
element $a\in\underline{g}:=\mathrm{Lie}(G)$, defines a Lie
derivative $L_{\xi^{a}}$ on $\Gamma (M,L)$ by the formula
$L_{\xi^{a}}(s):=\frac{d}{dt}|_{t=0}\mathrm{exp}(-ta)\cdot s$. The
universal contact form $\eta :TM\rightarrow L$ commutes with the
actions of $G$ on $M$ and $L$ and the $L$-valued $1$-form
$L_{\xi^{a}}(\eta)$ on $M$, defined by $L_{\xi^{a}}(\eta
)(X):=L_{\xi^{a}}(\eta (X))-\eta
(L_{\xi^{a}}X)$ for $X\in{\mathcal X}(M)$, is identically zero.\\

We briefly recall the reduction of CR manifolds as developed in
\cite{loose1}, \cite{loose2}. Let $G$ act by CR automorphisms on
the CR manifold $(M,H,I)$, $\Theta :M\rightarrow
L\otimes\underline{g}^{*}$ (where $\underline{g}^{*}$ denotes the
trivial bundle over $M$ with fiber
$\underline{g}^{*}:=\mathrm{Lie}(G)^{*}$) be the moment map,
defined by $\Theta_{a}:=\eta (\xi^{a})$ for every
$a\in\underline{g}$ and let $S:=\Theta^{-1}(\underline{0})$ (here
and elsewhere $\underline{0}$ will denote the image of the zero
section of a vector bundle). The Lie group $G$ preserves $S$ and
we can consider $p:S\rightarrow \hat{M}:=S/G$ the natural
projection. If the action of $G$ on $S$ is free and proper,
$\hat{M}$ is a CR manifold with the contact bundle
$\hat{H}:=p_{*}(TS\cap H)$. If $\theta$ is an arbitrary contact
form on $M$ and $g:=\frac{1}{2}d\theta (\cdot ,I\cdot )$ is the
corresponding metric on $H$, then the vector space $H_{x}$ (for
$x\in S$) has the $g$-orthogonal decomposition (independent of the
chosen contact form $\theta$)
\begin{equation}\label{descomps}
H_{x}=\mathcal T_{x}\oplus I\mathcal T_{x}\oplus E_{x},
\end{equation}
where ${\mathcal T}_{x}$ is the tangent space to the orbit of $G$ through $x$.
The space $E_{x}$ is $I$-invariant, is isomorphic via the
differential $p_{*}$ of $p$ with $\hat{H}_{p(x)}$ and its complex
structure induces the complex structure $\hat{I}$ of
$\hat{H}_{p(x)}.$\\

In this Section we add an horizontal and $G$-invariant Weyl connection
to this construction (see Definition below) and we study how it behaves 
under reduction. The Lie group $G$ is not necessarily connected.

\begin{defn}\label{condweyl}
Let $G$ be a Lie group acting by CR automorphisms on the CR-Weyl
manifold $(M,H,I,D)$. The Lie group $G$ acts by CR-Weyl
automorphisms if the following two conditions hold:

\begin{enumerate}

\item The connection $D$ is horizontal, i.e. it satisfies
\begin{equation}\label{inf}
D_{\xi^{a}}(s)=L_{\xi^{a}}(s),
\end{equation}
for every $s\in\Gamma (M,L)$ and vector field $\xi^{a}$ on $M$
generated by $G$.

\item The connection $D$ is $G$-invariant, i.e.
it satisfies
\begin{equation}\label{invariante}
g\cdot (D_{X}s)=D_{g_{*}X}(g\cdot{s}),
\end{equation}
for $g\in G$, $X\in{\mathcal X}(M)$, $s\in\Gamma (M,L)$.
\end{enumerate}
\end{defn}

\paragraph{Remark:}
Before studying the reduction, note that if a Weyl connection $D$ is 
horizontal and $G$ is connected, 
then $D$ is $G$-invariant if and only if its Faraday curvature $F^{D}$ 
satisfies
$F^{D}(\xi^{a},X)=0$, for every $X\in{\mathcal X}(M)$ and vector
field $\xi^{a}$ on $M$ generated by $G$: indeed, relation (\ref{inf})
implies
\begin{align*}
F^{D}(\xi^{a},X)&=L_{\xi^{a}}D_{X}(s)-D_{X}L_{\xi^{a}}(s)-D_{[\xi^{a},X]}(s)\\
&=\frac{d}{dt}\vert_{t=0}\mathrm{exp}(-ta)\cdot
D_{\mathrm{exp}(ta)_{*}(X)}\left(\mathrm{exp}(ta)\cdot s\right) ,
\end{align*}
from where we deduce our claim.\\

The following Proposition gives a nice description of the tangent
space to $S$ in terms of the Weyl connection $D$ (see also
\cite{loose1} for an equivalent proof).

\begin{prop}\label{tangenttoS}
For every $x\in S$, $T_{x}S=\{X\in T_{x}M|\quad d^{D}\eta
(X,\xi^{a})=0,\forall a\in\underline{g}\} .$
\end{prop}

\begin{proof}
A vector field $X$ is tangent to $S$ along $S$ if, for every
$a\in\underline{g}$, $(\Theta_{a})_{*}(X)$ is tangent to the image
of the zero section of $L$. Equivalently, if
$v^{D}(\Theta_{a})_{*}(X)=D_{X}(\Theta_{a})$ is zero, where we
recall that $v^{D}$ is the projection onto the last component in
the decomposition (\ref{desc}). On the other hand,
\begin{align*}
(d^{D}\eta )(X,\xi^{a})&=D_{X}(\Theta_{a})-D_{\xi^{a}}\left(\eta
(X)\right)
-\eta ([X,\xi^{a}])\\
&=D_{X}(\Theta_{a})-L_{\xi^{a}}(\eta )(X)=D_{X}(\Theta_{a}),\\
\end{align*}
where we have used $L_{\xi^{a}}(\eta )=0$ and the horizontality of
$D$.
\end{proof}

\begin{prop}\label{conweyl}
The connection $D$ induces a Weyl
connection $\hat{D}$ on the CR reduction
$(\hat{M},\hat{H},\hat{I})$. If $D$ is closed, so is $\hat{D}.$
\end{prop}

\begin{proof}
By dimension reasons, for every $x\in S$, the map
$$
\tilde{L}_{x}:=T_{x}(S)/T_{x}(S)\cap H_{x}\rightarrow L_{x}
$$
induced by the inclusion of $T_{x}(S)$ into $T_{x}(M)$ is an
isomorphism. Let $\tilde{L}\rightarrow S$ be the bundle whose
fiber over $x\in S$ is $\tilde{L}_{x}$. We notice that $\tilde{L}$
has a connection $\tilde{D}$ induced by the Weyl connection $D$:
if $i:S\rightarrow M$ is the inclusion, $\tilde{D}$ is the
connection $i^{*}D$ on $i^{*}L\cong\tilde{L}.$ The Lie group $G$
preserves $S$ and $H$ and induces an action on $\tilde{L}.$ The
differential of the natural projection $p:S\rightarrow \hat{M}$
induces bijective maps
$p_{*}:\tilde{L}_{x}\rightarrow\hat{L}_{p(x)}=T_{p(x)}\hat{M}/\hat{H}_{p(x)}$
and it is easy to see that the level sets of
$p_{*}:\tilde{L}\rightarrow\hat{L}$ are the orbits of $G$ on
$\tilde{L}$. It follows that the sections of the co-contact bundle
$\hat{L}$ of $\hat{M}$ can be identified with the $G$-invariant
sections of $\tilde{L}.$ Consider now a section $\hat{s}$ of
$\hat{L}$ and $\tilde{s}$ the (unique, $G$-invariant) section of
$\tilde{L}$ which projects to $\hat{s}$. Let $\hat{X}$ be a vector
field on $\hat{M}$ and $X_{1}$, $X_{2}$ two lifts of $X$ in $TS$.
Then $X_{1}-X_{2}$ is tangent to the orbits of $G$ in $S$ and the
horizontality of $D$ implies that
$\tilde{D}_{X_{1}}(\tilde{s})=\tilde{D}_{X_{2}}(\tilde{s}).$ Since
$\tilde{D}$, $X_{i}$ and $\tilde{s}$ are $G$-invariant, so is
$\tilde{D}_{X_{i}}(\tilde{s})$. The section of $\hat{L}$ induced by
$\tilde{D}_{X_{i}}(s)$ will be, by definition,
$\hat{D}_{\hat{X}}(\hat{s}).$ Since $p^{*}\hat{D}=\tilde{D}$,
$\hat{D}$ is closed if $D$ is.
\end{proof}

\begin{lem}\label{reeb}
The Reeb vector field $\psi^{D}$ is tangent to $S$ along $S$ and
projects to the Reeb vector field $\psi^{\hat{D}}$ of $\hat{D}.$
\end{lem}

\begin{proof}
Proposition \ref{tangenttoS} implies that $\psi^{D}\in
\tilde{L}^{-1}\otimes TS$ along $S$, because
$i_{\psi^{D}}d^{D}\eta=0$. Moreover, since $\psi^{D}$ is
$G$-invariant (because $\eta :TM\rightarrow L$ and $D$ are) with
respect to the action of $G$ on $\tilde{L}^{-1}\otimes TS$ , it is
also projectable onto $\hat{L}^{-1}\otimes T\hat{M}.$ Using the
natural isomorphism between $p^{*}\hat{L}$ and ${L}\vert_{S}$ it
is obvious that $p^{*}\hat{\eta}=\eta |_{TS}$ and
$p^{*}d^{\hat{D}}\hat{\eta}=d^{D}\eta |_{TS\times TS}.$ The
conclusion follows.
\end{proof}

\begin{thm}\label{sasakired}
Let $G$ be a Lie group acting by CR-Weyl automorphisms on a
Sasaki--Weyl manifold $(M,H,I,D)$. Let $\Theta :M\rightarrow
L\otimes\underline{g}^{*}$ defined by $\Theta_{a}:=\eta (\xi^{a})$
(for $a\in\underline{g}$) be the moment map. Suppose that $G$ acts
freely and properly on $S:=\Theta^{-1}(\underline{0})$ and let
$(\hat{M}:=S/G,\hat{H},\hat{I})$ be the CR quotient. Then
$M//G:=(\hat{M},\hat{H},\hat{I},\hat{D})$ is Sasaki--Weyl.
\end{thm}

\begin{proof}

Let $\hat{X}$ be a section of $\hat{H}$ and $X$ its unique lift on
$S$ which belongs to $E$ (recall the decomposition
(\ref{descomps}) of $H$ along $S$). We first claim that
$L_{\psi^{D}}(X)$ is a section of $L^{-1}\otimes E$. Indeed,
$L_{\psi^{D}}(X)$ belongs to $L^{-1}\otimes H$ because $X$ belongs
to $H$ and
\begin{align*}
(d^{D}\eta )\left( IL_{\psi^{D}}(X),\xi^{a}\right)&= (d^{D}\eta
)\left( L_{\psi^{D}}(IX),\xi^{a}\right) =-L_{\psi^{D}}(d^{D}\eta ) 
(IX,\xi^{a})\\
&=-i_{\psi^{D}}(F^{D}\wedge\eta )(IX,\xi^{a}) =-F^{D}(IX,\xi^{a})
\end{align*}
which is zero since $D$ is horizontal and $G$-invariant. (In the
above calculation we have used the Sasaki--Weyl condition on $D$,
$L_{\psi^{D}}(\xi^{a})=0$ because $\psi^{D}$ is $G$-invariant,
$(d^{D}\eta )(IX,\xi^{a})=0$ because $X$ belongs to $E$, and $\eta
(\xi^{a})=0$ along $S$). On the other hand, since $\psi^{D}$ and
$X$ are tangent to $S$, $L_{\psi^{D}}(X)$ is also tangent to $S$
and then, from Proposition \ref{tangenttoS}, $(d^{D}\eta
)(L_{\psi^{D}}(X),\xi^{a})=0.$ Our claim follows. It is easy to
see now that
$L_{\psi^{\hat{D}}}(\hat{I})(\hat{X})=p_{*}{L_{\psi^{D}}(I)(X)}$
(where $p_{*}:\tilde{L}^{-1}\otimes TS\rightarrow
(\tilde{L}^{-1}\otimes TS)/G=\hat{L}^{-1}\otimes T\hat{M}$ is the
natural projection) and hence $L_{\psi^{\hat{D}}}(\hat{I})=0$
because $D$ is Sasaki--Weyl.
\end{proof}

While $\hat{D}$ is closed when $D$ is, there are situations when
$D$ is exact but $\hat{D}$ is not. This is illustrated by
Proposition \ref{uni}, which is followed by a criteria (see
Proposition \ref{exact-non}) which expresses when the exactness
property of Sasaki--Weyl connections is preserved through the
reduction process.

\begin{prop}\label{uni}
Let $(M,H,I,D)$ be a Sasaki--Weyl manifold with a closed (not
necessarily exact) Weyl connection $D$. Let $\tilde{M}$ be the
universal cover of $M$. Then $\tilde{M}$ has an induced
Sasaki--Weyl structure $(\tilde{H},\tilde{I},\tilde{D}).$ Moreover,
$\tilde{D}$ is exact and $(M,H,I,D)$ is the Sasaki--Weyl reduction
of $(\tilde{M},\tilde{H},\tilde{I},\tilde{D})$ under the action of
$\pi_{1}(M)$ on $\tilde{M}.$
\end{prop}

\begin{proof}
Let $p:\tilde{M}\rightarrow M$ be the universal cover
with the Deck group $\Gamma$, isomorphic with $\pi_{1}(M).$ The CR
structure $(\tilde{H},\tilde{I})$ of $\tilde{M}$ is defined such
that the map $p:(\tilde{M},\tilde{H},\tilde{I})\rightarrow
(M,H,I)$ is a CR map. Note that the co-contact bundle $\tilde{L}$
of $(\tilde{M},\tilde{H},\tilde{I})$ is isomorphic with $p^{*}L$
and has the connection $\tilde{D}:=p^{*}D.$ It is clear that
$(\tilde{M},\tilde{H},\tilde{I},\tilde{D})$ is Sasaki--Weyl 
(the map $p$ is a local diffeomorphism of CR manifolds which preserves
the Weyl connections) and
that $\Gamma$ acts by CR automorphisms on $\tilde{M}.$ Moreover
$\tilde{D}$ is $\Gamma$-invariant: if $g\in\Gamma$, $s\in\Gamma
(M,L)$ and $X\in{\mathcal X}(M)$ has the lift
$\tilde{X}\in{\mathcal X}(\tilde{M})$, then
$$
g\cdot\tilde{D}_{\tilde{X}}(p^{*}s)=g\cdot p^{*}{D}_{X}(s)=p^{*}
D_{X}(s)=\tilde{D}_{\tilde{X}}(p^{*}s).
$$
The horizontality of $\tilde{D}$ with respect to the
$\Gamma$-action on $\tilde{M}$ is trivially satisfied, since
$\Gamma$ is discrete. It follows that $\Gamma$ acts by CR-Weyl
automorphisms on $(\tilde{M},\tilde{I},\tilde{H},\tilde{D})$ and
that the corresponding Sasaki--Weyl reduction is $(M,H,I,D)$ (the
moment map being trivial).
\end{proof}

Consider again the situation of Theorem \ref{sasakired}. Since $D$
is $G$ invariant and $M$ is connected, the action of $g\in G$ on
any global $D$-parallel (non-trivial) section $\mu$ of $L$ is of
the form $g\cdot \mu =\rho (g)\mu $, where $\rho (g)\in\mathbb{R}$
is independent of the chosen $\mu$. We obtain in this way a group
homomorphism $\rho :G\rightarrow \mathbb{R}\setminus\{ 0\}$
associated to the action of $G$ on $M$.

\begin{prop}\label{exact-non}
Suppose that the Weyl connection $D$ is exact and that the Lie group
$G$ is connected. Then the Weyl
connection $\hat{D}$ is also exact if and only if $\rho (G)=\{
1\}.$
\end{prop}

\begin{proof}
We shall use the notations employed in the proof of Proposition
\ref{conweyl}. Suppose first that $\rho (G)=\{ 1\}$ and take a
$D$-parallel (non-trivial) section $\mu$ of $L$. Then its
restriction to $S$ is $\tilde{D}$-parallel, $G$-invariant and it
determines a non-vanishing section $\hat{\mu}$ on $\hat{L}$,
which, from the definition of $\hat{D}$, must be
$\hat{D}$-parallel. For the converse, suppose now that $\hat{D}$
is exact, choose a $\hat{D}$-parallel non-vanishing section
$\hat{\mu}$ of $\hat{L}$, and let $\tilde{\mu}$ be its
corresponding $\tilde{D}$-parallel section of $\tilde{L}.$ Let
$x\in S$, ${\mathcal O}_{x}\subset S$ the orbit of $G$ through $x$
and $g\in G.$ Since $D$ is exact, there is a $D$-parallel section
$\mu$ of $L$ such that $\mu (x)=\tilde{\mu}(x)$. In particular,
the restriction of $\mu$ to $S$ is $\tilde{D}$-parallel, $\mu
=\tilde{\mu}$ on ${\mathcal O}_{x}$ (because ${\mathcal O}_{x}$ is
connected) and $g\cdot\mu =\mu$ on ${\mathcal O}_{x}$ since
$\tilde{\mu}$ is $G$-invariant. On the other hand, since $G$
preserves $D$, $g\cdot\mu$ is also $D$-parallel and it must
coincide with $\mu$ on the whole of $M$ (which is connected). This
means that $g\cdot \mu =\mu$, or $\rho (g)=1.$
\end{proof}

\paragraph{Remark:}With the notations of Proposition
\ref{exact-non}, notice that $\theta :=\mu^{-1}\eta$ is a contact
form of a Sasaki structure and is also $G$-invariant. A Sasaki
reduction of the Sasaki manifold $(M,H,I,\theta )$ is therefore
defined (see \cite{orneag}), whose underlying CR manifold is
$(\hat{M},\hat{H},\hat{I})$ and whose contact form is, as it can
be readily checked, the form $\hat{\mu}^{-1}\hat{\eta}.$\\

We end this Section by stating the commutativity between the
Sasaki--Weyl reductions and the operation of taking universal
covers. This is expressed by the following Proposition, whose
proof, lengthy but straightforward, will be omitted. Recall first
that an action of a Lie group $G$ on a manifold $M$ always lifts 
to an action of the universal cover $\tilde{G}$ of $G$ on the
universal cover $\tilde{M}$ of $M$, which commutes with the action
of $\pi_{1}(M)$ on $\tilde{M}.$

\begin{prop}
Let $(M,H,I,D)$ be a Sasaki--Weyl manifold and $G$ a Lie group
which acts by CR-Weyl automorphisms on $M$. Let $\tilde{G}$ be the
universal cover of $G$. The following statements hold:
\begin{enumerate}
\item The group $\tilde{G}$ acts by CR-Weyl automorphisms on
the universal cover $(\tilde{M},\tilde{H},\tilde{I},\tilde{D})$ of
$M$.

\item The action of $\pi_{1}(M)$ on $\tilde{M}$ induces a CR-Weyl action of
$\pi_{1}({M})$ on $\tilde{M}//\tilde{G}$, and there is a CR
diffeomorphism $\tilde{M}//\tilde{G}/\pi_{1}(M)\cong M//G$ which
preserves the Weyl connections.
\end{enumerate}

\end{prop}

\section{The commutativity with the cone
construction}\label{comutativitate}

In this Section we consider a connected Lie group $G$ acting by orientation
preserving CR-Weyl automorphisms on the Sasaki--Weyl manifold
$(M,H,I,D)$ whose cone $(L^{+},\omega ,J)$ is l.c.K. In
particular, the curvature $F^{D}$ of $D$ must satisfy
$i_{\xi^{a}}F^{D}=0$ for every vector field $\xi^{a}$ on $M$
generated by $G$, and $F^{D}=kd^{D}\eta$, for a section $k$ of
$L^{-1}$. Note that these two conditions on the curvature $F^{D}$ actually 
imply that $F^{D}$ is identically zero on the complement of the zero set
$S:=\Theta^{-1}(\underline{0})$ of the moment map $\Theta$ of the action
of $G$ on $(M,H,I)$. Indeed, suppose, on the contrary, that $k$ is non-zero 
at a point $x\in M\setminus S$
and let $a\in\underline{g}$ such that $\eta (\xi^{a})$ is non-zero on a small
neighborhood $V$ of $x$, included in $M\setminus S$. From the proof of
Proposition \ref{tangenttoS} we know that 
$$
F^{D}(\xi^{a},X)=kD_{X}\left(\eta (\xi^{a})\right) ,\quad\forall X\in 
{\mathcal X}(M). 
$$
Since $F^{D}(X,\xi^{a})=0$ we deduce that $\eta (\xi^{a})$ is $D$-parallel 
and non-vanishing
on $V$. It follows that $F^{D}$ is zero on $V$, and hence $k(x)=0$. 
We have obtained a contradiction.\\

For the rest of this section we will suppose that $D$ is closed. In particular,
the canonical Weyl connection of the l.c.K. cone
$(L^{+},\omega ,J)$ is $\bar{D}:=\pi^{*}D$. Since the action of $G$ on $M$ 
is orientation preserving, it
induces an action on $L^{+}$ and on the bundle $\pi^{*}L$ on
$L^{+}.$ We shall denote by $\hat{\xi}^{a}$ the vector field on
$L^{+}$ generated by $a\in\underline{g}$.

\begin{lem}\label{liftoriz}
The vector fields generated by the action of $G$ on $L^{+}$ are
$D$-horizontal lifts of the vector fields generated by the action
of $G$ on $M$.
\end{lem}

\begin{proof}
Let $s\in\Gamma (M,L)$ and define $s_{t}:=\mathrm{exp}(-ta)\cdot
s$. Since $D$ is horizontal, we infer that
$\frac{d}{dt}|_{t=0}s_{t}=D_{\xi^{a}}(s).$ Our claim readily
follows: for $x\in M$,
\begin{align*}
\widetilde{{\xi}^{a}}_{s (x)}&=s_{*}({\xi}^{a}_{x})-
D_{{\xi}^{a}_{x}}(s )\\
&=s_{*}\left(\frac{d}{dt}|_{t=0}\mathrm{exp}(ta) x\right)
-D_{{\xi}^{a}_{x}}(s)\\
&=\frac{d}{dt}|_{t=0}{\mathrm{exp}(ta)}\cdot
s_{t}(x) -D_{\xi^{a}_{x}}(s)\\
&=\widehat{\xi^{a}}_{s(x)}+\frac{d}{dt}|_{t=0}s_{t}(x)-D_{\xi^{a}_{x}}(s)\\
&=\widehat{\xi^{a}}_{s(x)}.
\end{align*}

\end{proof}

\begin{lem}
The connection $\bar{D}$ on the bundle $\pi^{*}L\rightarrow L^{+}$
is $G$-invariant and horizontal.
\end{lem}

\begin{proof}
Since $\bar{D}$ is closed and $G$ is connected, it is enough to
show that $\bar{D}$ is horizontal. The very definition of the
action of $G$ on $\pi^{*}L$ implies
$\mathrm{exp}(-ta)\cdot\pi^{*}(s)=\pi^{*}(\mathrm{exp}(-ta)\cdot
s)$, for $s\in\Gamma (M,L)$ and $a\in\underline{g}.$ It is enough
to take the derivative (with respect to $t$) of this relation, to
use the horizontality of $D$ and Lemma \ref{liftoriz}.
\end{proof}

\begin{prop}\label{existence}
The Lie group $G$ preserves the hermitian conformal structure of
$(L^{+},\omega ,J)$
\end{prop}

\begin{proof}
We first prove that for every $X\in{\mathcal X}(M)$ and $g\in G$,
$g_{*}\tilde{X}=\widetilde{g_{*}X}$: for this, let $s\in\Gamma
(M,L)$ and $x\in M.$ The $G$-invariance of $D$ implies
$$
g_{*}\tilde{X}_{s(x)}=g_{*}\left( s_{*}X_{x}-D_{X_{x}}s\right)
=(g\cdot s)_{*}g_{*}(X_{x})-D_{g_{*}X_{x}}(g\cdot
s)=\widetilde{g_{*}(X_{x})}_{g\cdot s(x)}
$$
which is our claim. It is now easy to see that $G$ acts by
holomorphic transformations on $L^{+}:$ for every $X\in H$,
$$
g_{*}(J\tilde{X})=g_{*}(\widetilde{IX})=\widetilde{g_{*}(IX)}=
\widetilde{Ig_{*}(X)}=Jg_{*}(\tilde{X}).
$$
In a similar way, $g_{*}J(\psi^{D}(s_{0}))=g_{*}s_{0}= J g_{*}\left(
\psi^{D}(s_{0})\right)$ because $\psi^{D}$ is $G$-invariant. It
follows that $g_{*}\left( J\psi^{D}\right) =J\left(
g_{*}\psi^{D}\right)$ and that $G$ acts by bi-holomorphic
transformations on $L^{+}.$ Since $\bar{D}$ and
$\sigma\pi^{*}\eta$ are $G$-invariant, it is clear that $\omega
=\frac{1}{2}d^{\bar{D}}(\sigma\pi^{*}\eta )$ is also
$G$-invariant: $g\cdot\omega (X,Y)=\omega (g_{*}(X),g_{*}(Y))$,
for every $g\in G$ and $X,Y\in{\mathcal X}(L^{+}).$ The conclusion
follows.

\end{proof}

According to Proposition \ref{existence}, there is a l.c.K.
reduction (see \cite{gaud} or \cite{partonornea}) of
$(L^{+},\omega ,J)$ under the action of $G$. It has a
distinguished moment map provided by the following Proposition.

\begin{prop}\label{moment}
Let ${\Theta}:M\rightarrow L\otimes\underline{g}^{*}$ be the
moment map for the action of $G$ on $(M,H,I,D)$. Then the map
$\tilde{\Theta}
:L^{+}\rightarrow\pi^{*}(L^{2}\otimes\underline{g}^{*})$ defined
by $\tilde{\Theta}_{a}:=\frac{1}{2}\pi^{*}({\Theta}_{a})\sigma$
for $a\in\underline{g}$, is a moment map for the action on $G$ on
the l.c.K. manifold $(L^{+},J, \omega ).$
\end{prop}

\begin{proof}
Let $X\in {\mathcal X}(L^{+})$ be $\pi$-projectable and
$Y:=\pi_{*}(X).$ From $\omega
=\frac{1}{2}d^{\bar{D}}(\sigma\pi^{*}\eta )$, we get
\begin{align*}
\omega (X,\tilde{\xi}^{a})&=\frac{1}{2}\bar{D}_{X}\left(
\sigma\pi^{*}\Theta_{a}\right)
-\frac{1}{2}\bar{D}_{\tilde{\xi}^{a}}\left(\sigma\pi^{*}\eta
(Y)\right)-\frac{\sigma}{2}\pi^{*}\eta ([Y,\xi^{a}])\\
&=\frac{1}{2}\bar{D}_{X}\left(\sigma\pi^{*}\Theta_{a}\right)
-\frac{1}{2}\bar{D}_{\tilde{\xi}^{a}}(\sigma )\pi^{*}\eta
(Y)-\frac{\sigma}{2} \bar{D}_{\tilde{\xi}^{a}}\left(\pi^{*}\eta
(Y)\right)-\frac{\sigma}{2}\pi^{*}\eta ([Y,\xi^{a}]).
\end{align*}
On the other hand, since $D$ is horizontal and $\eta$ is
$G$-invariant,
$$
\bar{D}_{\tilde{\xi}^{a}}\left(\pi^{*}\eta (Y)\right) =
\pi^{*}L_{\xi^{a}}\left(\eta (Y)\right) =\pi^{*}\eta
([\xi^{a},Y]).
$$
Also, since $\bar{D}$ is horizontal and $\sigma$ is $G$-invariant,
$\bar{D}_{\tilde{\xi}^{a}}(\sigma )=L_{\tilde{\xi}^{a}}(\sigma
)=0.$ We have proved that $\omega
(X,\tilde{\xi}^{a})=\frac{1}{2}\bar{D}_{X}\left(\sigma\pi^{*}\Theta_{a}\right)$
which implies the conclusion.
\end{proof}

\begin{thm}
Let $(M,H,I,D)$ be a Sasaki--Weyl manifold whose Weyl connection 
$D$ is closed. Let $G$
be a connected Lie group acting by orientation preserving CR-Weyl
automorphisms on $M$ such that the underlying manifold of the
Sasaki--Weyl reduction $M//G=(\hat{M},\hat{H},\hat{I},\hat{D})$ is
smooth. Then the l.c.K. reduction of $(L^{+},\omega ,J)$ under the
action of $G$ coincides with l.c.K. cone of the Sasaki--Weyl
reduction $M//G$.
\end{thm}

\begin{proof}
We shall use the following maps: $\pi :L\rightarrow
M$, $\hat{\pi}:\hat{L}\rightarrow\hat{M}$ for the co-contact bundles
of $M$ and $\hat{M}$ respectively,
$\hat{p}:\tilde{L}^{+}\rightarrow\hat{L}^{+}$,
$p:S\rightarrow\hat{M}$ the natural projections, where $\tilde{L}$
is the restriction of the bundle $L$ to the zero set
$S:=\Theta^{-1}(\underline{0})$ of the moment map $\Theta$ of the
action of $G$ on $M.$ Proposition \ref{moment} readily implies
that $\tilde{L}^{+}$ is the zero set of the moment map
$\tilde{\Theta}$ of $G$ on $L^{+}$  and that the underlying manifold of 
the l.c.K.
reduction of $(L^{+},\omega ,J)$ is $\tilde{L}^{+}/G$ and
coincides with the cone $\hat{L}^{+}$ of $(\hat{M},\hat{H},\hat{I},\hat{D}).$ 
We will now
check that the complex structure $J_{1}$ of $\hat{L}^{+}$, induced
by the l.c.K. reduction, coincides with the cone complex structure
$J_{2}$ of $\hat{L}^{+}$. For every $y\in\tilde{L}^{+}$ there is a
decomposition
$$
T_{y}L^{+}=L_{\pi (y)}\oplus\widetilde{\mathcal T_{\pi (y)}}^{D}\oplus
\widetilde {I{\mathcal T}_{\pi (y)}}^{D}\oplus\widetilde{E_{\pi
(y)}}^{D}\oplus \widetilde{\Lambda_{\pi (y)}}^{D},
$$
where the superscript \lq\lq$D$'' denotes horizontal lift with respect
to the connection $D$ at the point $y$, $\mathcal T_{\pi (y)}$ is the 
tangent space to the
orbit of $G$ through the point $\pi (y)$ of $M$ and $E_{\pi (y)}$
is the hermitian orthogonal complement of $\mathcal T_{\pi (y)}\oplus
I{\mathcal T_{\pi (y)}}$ in $H_{\pi (y)}$ (with respect to the metric of
$H_{\pi (y)}$ associated to any contact form of the contact bundle $H$). Lemma
\ref{liftoriz} implies that $\widetilde{\mathcal T_{\pi (y)}}^{D}$ is the
tangent space to the orbit of $G$ in $L^{+}$ through $y$ and the
definition of the cone complex structure $J$ of $L^{+}$ implies that
$\widetilde {I{\mathcal T}_{\pi (y)}}^{D}=J 
\widetilde{\mathcal T_{\pi (y)}}^{D}.$ The
differential $\hat{p}_{*}:L_{\pi (y)}\oplus \widetilde{E_{\pi
(y)}}^{D}\oplus \widetilde{\Lambda_{\pi (y)}}^{D}
\rightarrow T_{\hat{p}(y)}\hat{L}^{+}$ of $\hat{p}$ is an
isomorphism and the complex structure ${J}_{1}$ of the l.c.K reduction is
\begin{align*}
{J}_{1}\left( \hat{p}_{*}\widetilde{\psi^{D}_{\pi
(y)}(y)}^{D}\right) &=-\hat{p}(y)\in T^{V}_{\hat{p}(y)}\hat{L}^{+}\\
{J}_{1}\left(\hat{p}_{*}\widetilde{X}^{D}\right)&=\hat{p}_{*}\widetilde{IX}^{D},\quad
\forall X\in E_{\pi (y)}.
\end{align*}
On the other hand, since $p_{*}E_{\pi (y)}=\hat{H}_{p\pi
(y)}=\hat{H}_{\hat{\pi}\hat{p}(y)}$, the Weyl connection $\hat{D}$
induces the decomposition
$$
T_{\hat{p}(y)}\hat{L}=\hat{L}_{\hat{\pi}\hat{p}(y)} \oplus
\widetilde{\hat{\Lambda}_{\hat{\pi}\hat{p}
(y)}}^{\hat{D}} \oplus\widetilde{p_{*}E_{\pi (y) }}^{\hat{D}},
$$
where the superscript \lq\lq$\hat{D}$'' denotes $\hat{D}$-horizontal
lift and $\hat{\Lambda}$ is the image of the Reeb vector field
$\psi^{\hat{D}}:\hat{L}\rightarrow T\hat{M}$. Then 
the cone complex structure $J_{2}$ of
$\hat{L}$ can be written down explicitely as follows:
\begin{align*}
{J}_{2} \left(\widetilde{\psi^{\hat{D}}_{\hat{\pi}(\hat{p}(y))}\hat{p}(y)}^{\hat{D}}\right) &=-\hat{p}(y)\in T^{V}_{\hat{p}(y)}(\hat{L}^{+})\\
{J}_{2}(\tilde{Y}^{\hat{D}})&=\widetilde{\hat{I}Y}^{\hat{D}}
,\quad \forall Y\in p_{*}E_{\pi (y)}.
\end{align*}
In order to show that $J_{1}=J_{2}$, it is enough to prove that
\begin{equation}\label{coincidence}
\hat{p}_{*}(\tilde{X}^{D}_{y})
=(\widetilde{p_{*}X}^{\hat{D}})_{\hat{p}(y)},\quad\forall  X\in
T_{\pi (y)}S,
\end{equation}
(recall that $p_{*}\psi^{D}=\psi^{\hat{D}}$, from Lemma
\ref{reeb}), which can be done as follows: let $\hat{s}\in\Gamma
(\hat{M},\hat{L})$ such that
$\hat{s}\hat{\pi}\hat{p}(y)=\hat{p}(y)$ and consider the
corresponding ($G$-invariant) section $s$ of $\tilde{L}.$ Then
$\hat{p}\circ s =\hat{s}\circ p$, $s\pi (y)=y$ and
\begin{align*}
\hat{p}_{*}\left(\tilde{X}^{D}_{y}\right)&=\hat{p}_{*}\left(s_{*}(X_{\pi (y)})-D_{X_{\pi (y)}}(s)\right)\\
&=\hat{s}_{*}p_{*}(X_{\pi (y)})-\hat{D}_{p_{*}(X_{\pi (y)})}(\hat{s})\\
&=\widetilde{p_{*}(X_{\pi (y)})}^{\hat{D}}_{\hat{p}(y)},
\end{align*}
which implies (\ref{coincidence}). The coincidence of the two
metrics on $\hat{L}^{+}$ can be proved in a similar way (note that
the metric coming from the l.c.K. reduction takes values in a
bundle ${\mathcal D}$ characterized by $\hat{p}^{*}\mathcal D\cong
\pi^{*}L^{2}|_{\tilde{L}^{+}}$; the cone metric of $\hat{L}^{+}$
takes values in $\hat{\pi}^{*}\hat{L}^{2}$; since
$\hat{p}^{*}\hat{\pi}^{*}\hat{L}^{2}=\pi^{*}p^{*}\hat{L}^{2}\cong\pi^{*}(L^{2}|_{S})
\cong\pi^{*}L^{2}|_{\tilde{L}^{+}}$, it follows that the two
metrics take values in the same line bundle).

\end{proof}

\section{Examples}

\begin{enumerate}

\item A class of non-closed Sasaki--Weyl connections can be
constructed in the following way: suppose that $D_{0}$ is an exact
Sasaki--Weyl connection on the CR manifold $(M,H,I)$ and let
$s_{0}$ be an arbitrary section of $L^{-1}.$ The connection
$D:=D_{0}+s_{0}\eta$ is non-closed and Sasaki--Weyl: $d^{D}\eta
=d^{D_{0}}\eta$ implies that $\psi^{D}=\psi^{D_{0}}$ and
$L_{\psi^{D}}(I)=L_{\psi^{D_{0}}}(I)=0$. It follows that any CR
manifold which admits a compatible Sasaki structure admits
infinitely many of non-closed Sasaki--Weyl connections,
parametrized by smooth sections of the dual of its co-contact
bundle. Suppose now that $s_{0}$ is $D_{0}$-parallel and
non-vanishing. Then $F^{D}=s_{0}d^{D}\eta$,
$D(s_{0})+s_{0}^{2}\eta =0$ and the cone of $(M,H,I,D)$ is
globally conformal K\"{a}hler (see Theorem \ref{lcK}), with the
K\"{a}hler structure defined explicitly as follows: let $\theta
:=s_{0}\eta$ and $T_{0}:=\psi^{D_{0}}(s_{0}^{-1})$ be the contact form
and Reeb vector field associated to $s_{0}^{-1}.$ Consider
the trivialization of $L^{+}$ determined by $s_{0}^{-1}$. The
$1$-form $\gamma$ from the proof of Lemma \ref{potential} is in this case
$s_{0}\eta$. It vanishes on $H$ and it takes the value one 
when applied to $T_{0}$. The complex structure $J$, 
written down explicitely in  the proof of Lemma \ref{potential}, becomes
$$
JX=IX;\quad J\left(\frac{\partial}{\partial t}\right) =T_{0}-
t\frac{\partial}{\partial t} . 
$$
The K\"{a}hler form is $\omega =\frac{1}{2}d(t\theta_{0}).$

\item We shall illustrate the Sasaki--Weyl reduction in an explicit
situation. Let $S^{1}\times\mathbb{Z}$ act on
$\bar{M}:=\mathbb{C}^{n}\setminus\{ 0\}\times\mathbb{R}^{>0}$ by
$$
(e^{i\theta}, m)\cdot (z_{1},\cdots ,z_{n},t):=
(e^{ia_{1}\theta}\lambda^{m}z_{1},\cdots ,
e^{ia_{n}\theta}\lambda^{m}z_{m},\lambda^{2m}t),
$$
where $\lambda >1$ and $a_{1},\cdots ,a_{n}$ are integers, not all
of the same sign. The
$1$-form $\theta :=\sum_{p=1}^{n}x_{p}dy_{p}-y_{p}dx_{p}-dt$ is
preserved, up to homothety, by the action of
$S^{1}\times\mathbb{Z}$ and determines a Sasaki structure on
$\bar{M}$.  The Sasaki--Weyl reduction of $\bar{M}$ (with this
Sasaki structure), under the action of $S^{1}\times\mathbb{Z}$ is
defined, and the zero set of its moment map is
$$
S:=\{ (z_{1},\cdots ,z_{n})\in\mathbb{C}^{n}\setminus\{ 0\}:
\sum_{p=1}^{n}a_{p}|z_{p}|^{2}=0\}\times\mathbb{R}^{>0}.
$$
It follows that the underlying manifold $\hat{M}_{\lambda
,a_{1},\cdots ,a_{n}}$ of the Sasaki--Weyl reduction is
diffeomorphic with $(S\cap S^{2n-1})/S^{1}\times(
\mathbb{R}^{>0}\times\mathbb{R}^{>0})/\mathbb{Z}$, where
$S^{2n-1}$ is the unit sphere in $\mathbb{C}^{n}$, $S^{1}$ acts on
$S^{2n-1}$ by
\begin{equation}\label{actiune}
e^{i\theta}\cdot (z_{1},\cdots ,z_{n}):=
(e^{ia_{1}\theta}z_{1},\cdots , e^{ia_{n}\theta}z_{n})
\end{equation}
and $\mathbb{Z}$ acts on $\mathbb{R}^{>0}\times\mathbb{R}^{>0}$ by
$m\cdot (t_{1},t_{2}):=(\lambda^{m}t_{1},\lambda^{2m}t_{2}).$ Note
that the space $M_{a_{1},\cdots ,a_{n}}:=(S\cap S^{2n-1})/S^{1}$
coincides with the Sasaki reduction of $S^{2n-1}$ with its natural
Sasaki structure, under the action (\ref{actiune}) of $S^{1}$, and
has been identified, for various $(a_{1},\cdots ,a_{n})$, in
\cite{orneag}. We obtain closed, non-exact Sasaki--Weyl connections
on products $M_{a_{1},\cdots ,a_{n}}\times S^{1}\times\mathbb{R}.$

\item The previous Example can be considerably generalized. Let
$(M,g,J)$ be a K\"{a}hler manifold with the K\"{a}hler form
$\omega =d\alpha$, for a $1$-form $\alpha$ on $M$. Let $G$ be a
discrete group acting freely, proper discontinuously, by
holomorphic transformations on $(M,g,J)$ and such that for every
$g\in G$, $g^{*}\alpha =\rho (g)\alpha$ for a positive number
$\rho (g).$ The $1$-form $\theta :=\alpha -dt$ determines a Sasaki
structure on $M\times\mathbb{R}^{>0}$ and is preserved, up to
homothety, by the action $g\cdot (x,t):=(g\cdot x,\rho (g) t)$ of
$G$ on $M\times\mathbb{R}^{>0}.$ Proposition \ref{exact-non}
implies that the quotient manifold $(M\times\mathbb{R}^{>0})/G$
has a closed non-exact Sasaki--Weyl connection.

Consider now a Sasaki manifold $(N,H,I,T)$ with Reeb vector field
$T$ and $(N\times\mathbb{R}^{>0},\omega ,J)$ its K\"{a}hler cone,
whose K\"{a}hler form has potential $t^{2}$ and whose complex
structure $J$ satisfies $J(V):=T$ (where
$V:=t\frac{\partial}{\partial t}$ is the radial vector field on
$\mathbb{R}^{>0}$) and $J|_{H}=I$ on $H$. The K\"{a}hler cone
$(N\times\mathbb{R}^{>0},\omega ,J)$ admits an holomorphic
$\mathbb{Z}$-action defined by
\begin{equation}\label{z}
m\cdot (x,t):=(x,\lambda^{m}t),
\end{equation}
for a fixed integer $\lambda >1$, which preserves $\alpha
:=\frac{1}{4}d^{J}t^{2}$ up to homothety. Therefore, we can take
$M:=(N\times\mathbb{R}^{>0},\omega ,J)$, $G:=\mathbb{Z}$ act on
$M$ by (\ref{z}), and deduce that the quotient $N\times
S^{1}\times S^{1}$ of $M\times\mathbb{R}^{>0}$ under the
$\mathbb{Z}$-action $m\cdot
(x,t_{1},t_{2}):=(x,\lambda^{m}t_{1},\lambda^{2m}t_{2})$ has a
closed, non-exact, Sasaki--Weyl connection. Moreover, if $K$ is a
Lie group which acts by CR automorphisms on $N$ and preserves $T$,
then we can couple the action of $K$ on $N$ with the
$\mathbb{Z}$-action (\ref{z}) on $\mathbb{R}^{>0}$, to get an
holomorphic action of $K\times\mathbb{Z}$ on
$N\times\mathbb{R}^{>0}$ which preserves, up to homothety, the
$1$-form $\alpha =\frac{1}{4}d^{J}t^{2}.$ Again, we can take
$M:=(N\times\mathbb{R}^{>0},\omega ,J)$ and
$G:=K\times\mathbb{Z}.$ A Sasaki--Weyl reduction
$N\times\mathbb{R}^{>0}\times\mathbb{R}^{>0}//K\times\mathbb{Z}$
is defined, which is isomorphic, as a Sasaki--Weyl manifold,  with
$N//K\times S^{1}\times\mathbb{R}$ (where $N//K$ is the Sasaki
reduction of $(N,H,I,T)$ under the action of $K$). Example 1 is
just a particular case of this construction: $K:=S^{1}$ acts on
$N:=S^{2n-1}$ by (\ref{actiune}) preserving its standard Sasaki
structure and $\mathbb{Z}$ acts on $\mathbb{R}^{>0}$ by (\ref{z}).

\end{enumerate}

\end{document}